# Langevin equations solved by Markov chains (instead of integrals), two solutions with multiplicative noise and their properties

Dietrich Ryter

ryterdm@gawnet.ch

Midartweg 3   CH-4500 Solothurn  Switzerland

*Abstract*

Langevin equations are usually understood as integral equations, since continuous and Markov solutions are nowhere differentiable. This is avoidable by the fact that white Gaussian noise can be seen as a generalized derivative of the Wiener process. The simplest LE is thus solved by Wiener paths, which can be computed by a Markov chain. That crucial property persists for functions of these paths, and solutions of the general LE are therefore assumed to have that form. In short time steps these functions must be linear or quadratic. Both cases lead to sampled versions of existing stochastic integrals ('Ito' and 'anti-Ito' respectively). In physics they apply for Brownian motion in presence of a temperature gradient: without a feedback to the medium in the linear case, and for a tagged molecule in the quadratic one; in the latter case the Fokker-Planck equation describes a local equilibrium.

## I. Introduction

Langevin equations are ordinary differential equations perturbed by noise. They arise in statistical physics and in other fields like e.g. in stochastic control. Since continuous and Markov solutions are nowhere differentiable, the LEs are usually understood as integral equations which are not necessarily solvable by increments in consecutive short time steps (this is crucial for numerical computations). The underlying 'stochastic integral' is not unique when the impact of the noise depends on the random state (i.e. for 'multiplicative noise' occurring in non-equilibrium systems): the defining Riemannian sum then depends on the evaluation point characterized by a free parameter $0 \leq \propto \leq 1$ ($\propto = 0$ for 'Ito', $\propto = 1/2$ for 'Stratonovich' and $\propto = 1$ for 'anti-Ito') [1-8]. The Markov property implies that $\propto = 0$ or 1. A further ambiguity arises when $\propto > 0$ : the discretization of time is then not unique since the stochastic calculus admits two different Riemannian sums converging to the same limit; one of them is restated by the new approach and - as it turns out - also singled out by the Fokker-Planck equation [9].

The present approach is quite different. It is based on the fact that the simplest LE is solved by the (vector) Wiener process [10], and it assumes that the resolving paths of the general LE are *functions* of the Wiener paths. Since the latter can be computed by a Markov chain, the same holds for any functions of them (evaluated at the same epochs). Short time steps admit either linear or quadratic functions. The linear case corresponds to the standard $\propto = 0$. For multiplicative noise the quadratic case is also possible (with $\propto = 1$). The analysis of the quadratic case is complicated by the fact that the increments are not Gaussian distributed. The corresponding FPE is discussed in the Appendix.

Physically, $\propto = 1$ applies for a 'tagged' molecule in a fluid with a temperature gradient, and the FPE indicates a local equilibrium with the other molecules.

## II. The Langevin equation, and first results

The Wiener process plays a central role: it is Gaussian distributed, with continuous but non-differentiable paths $W(t)$ obeying $< W(t) - W(0) > = 0$ and $< [W(t) - W(0)]^2 > = t$. It is Markovian, both for continuous and discrete time. A crucial finding [10] states that the 'derivative' of $W(t)$ can be identified with the (highly singular) white Gaussian noise $\xi(t)$. A deterministic dynamical system $\dot{\vec{x}} = \vec{a}(\vec{x})$ is thus easily supplied by noise, to become the Langevin equation

$$\dot{\vec{X}} = \vec{a}(\vec{X}) + \underline{B}(\vec{X})\,\vec{\xi} \quad \text{or} \quad \dot{X}^i = a^i(\vec{X}) + b^{ik}(\vec{X})\,\xi_k \tag{2.1}$$

for the resolving random paths $\vec{X}(t)$. The $\xi_k$ are independent. The vector Wiener process $\vec{W}(t)$ solves (2.1) when $\vec{a} = \vec{0}$ and $\underline{B} = \underline{I}$.

The remainder of the paper is focused on path increments in a short time $\tau$, and with a given initial value $\vec{X}(0) = \vec{x}$. For constant $\underline{B}$, and with $\Delta W_k := W_k(\tau) - W_k(0)$, they are linear functions of the $\Delta W_k$ :

$$\Delta X^i = a^i(\vec{x})\tau + b^{ik}\Delta W_k + o(\tau)\ . \tag{2.2}$$

A nonconstant $\underline{B}$ implies that the impact of the noise depends on the random $\vec{X}(t)$ ('multiplicative noise'). It is then *assumed* that $\Delta X^i$ is still a function of the $\Delta W_k$ – so that $\vec{X}(t)$ is determined by a Markov chain, at times separated by $\tau$. This allows a completion of (2.2) by products $\Delta W_i \Delta W_k$ (which are of the order $O(\tau)$), multiplied by derivatives of $b^{ik}$ at $\vec{x}$ (higher powers do not contribute since $< |\Delta W|^n >$ is proportional to $\tau^{n/2}$). Details will be given below.

The *covariance matrix* $<\Delta\vec{X}\,\Delta\vec{X}^T>$ plays an essential role. For its evaluation in the order $O(\tau)$ it is sufficient to consider the lowest approximation $\Delta\vec{X} = \underline{B}(\vec{x})\Delta\vec{W}$. Then the dyadic product $\Delta\vec{X}\,\Delta\vec{X}^T$ equals $\underline{B}(\vec{x})\,\Delta\vec{W}\Delta\vec{W}^T\underline{B}^T(\vec{x})$. By $< \Delta\vec{W}\Delta\vec{W}^T > = \underline{I}\tau$ the expectation has the required order $O(\tau)$ and is given by

$$< \Delta\vec{X}\,\Delta\vec{X}^T > = \underline{D}(\vec{x})\,\tau \qquad (2.3)$$

with the symmetric 'diffusion matrix' $\underline{D}(\vec{x}) \coloneqq \underline{B}(\vec{x})\underline{B}^T(\vec{x})$. Clearly the off-diagonal elements of $\underline{D}(\vec{x})$ describe the correlation between the components of $\Delta\vec{X}$. [Multiplying $\underline{B}(\vec{x})$ by a constant orthogonal $\underline{O}$ does not change $\underline{D}(\vec{x})$; it amounts to replacing $\vec{W}(t)$ by a stochastically equivalent Wiener process].

## III. Multiplicative noise

3.1 *The linear case*

The increment (2.2) with $b^{ik}$ taken at $\vec{x}$ satisfies the basic assumption:

$$\Delta X^i = a^i(\vec{x})\tau + b^{ik}(\vec{x})\,\Delta W_k + o(\tau)\,. \qquad (3.1)$$

It corresponds to the Ito case $\alpha = 0$, and it admits any smooth $\underline{B}(\vec{x})$. Clearly (3.1) is Gaussian distributed, with the mean $\vec{a}\,\tau$ and with the variance $\underline{D}\,\tau$ [there is an exception when $\vec{a}$ has an equilibrium point at the initial $\vec{x}$: then the local derivatives of $\vec{a}$ contribute to the variance, see [11,12] ].

The FPE for the density function $p(\vec{x},t)$ reads

$$p_{,t} = [-a^i p + (1/2)(D^{ik}p)_{,k}]_{,i}\ . \qquad (3.2)$$

It specifies a probability current $\vec{J}(\vec{x},t)$, given by

$$J^i = (a^i - D^{ik}{}_{,k}/2)\,p \quad \text{where} \quad \nabla p = \vec{0}\,, \qquad (3.3)$$

see (A.5) .

It is worthwhile to discuss (3.3) in a physical situation (in one dimension, and with $a = 0$, so that $J = (-D'/2)\,p$ where $p' = 0$ : Consider a Brownian particle (BP) in a medium with a temperature gradient ($D$ is proportional to the temperature). Then $J$ indicates a trend towards the cool side, due to elastic hits by the molecules of the substrate that are stronger and more frequent from the warm side. The BP must be massive (slow) enough that it does not change the energy of the molecules, otherwise a memory effect arises and excludes the Markov property.

3.2 *The quadratic case, local equilibrium*

Here $\vec{a}$ is omitted for simplicity, and will be included afterwards.

In one dimension the quadratic extension of (2.2) takes the form

$$\Delta X \approx b(x)\Delta W + b'(x)b(x)(\Delta W)^2 = \sqrt{D}\,\Delta W + (1/2)D'\,(\Delta W)^2 + o(\tau) \tag{3.4}$$

so that

$$< \Delta X > = (1/2)\,D'\,\tau + o(\tau)\ . \tag{3.5}$$

This corresponds to the anti-Ito case. In higher dimensions products $\Delta W_i \Delta W_k$ arise with $i \neq k$, they are of the order $O(\tau)$ and their mean vanishes. Their density function is nonsingular, in contrast to that with $i = k$ which diverges at zero. It follows that the $(\Delta W_i)^2$ predominate for small $\tau$, so that the $< \Delta X^i >$ evolve independently of each other; their growth is determined by the variation of $\underline{D}(\vec{x})$ in the respective direction. The off-diagonal elements of $\underline{D}(\vec{x})$ must be constant. This results in

$$\Delta X^i = b^{ik}\,\Delta W_k + (D^{ii}{}_{,i}/2)\,(\Delta W_i)^2 + o(\tau)$$

with the mean

$$< \Delta X^i > = (D^{ii}{}_{,i}/2)\,\tau + o(\tau)\ . \tag{3.6}$$

In the Appendix it is argued that a diagonal $\underline{D}$ entails a local equilibrium, which is crucial for the Markov property. The above increments then take the form

$$\Delta X^i = \sqrt{D^{ii}}\,\Delta W_i + (D^{ii}{}_{,i}/2)\,(\Delta W_i)^2 + o(\tau) \tag{3.7}$$

with the mean (3.6) . The impact of $r$-dimensional noise is thus specified by $r$ positive and independent functions $D^{ii}(x^i)$.

The density function of $\Delta\vec{X}$ is a product of one-dimensional ones, as given by (A.9) :

$$p(x,\tau) \approx [2\pi D(0)\tau]^{-1/2}\,[1 + x^3 D'(0) D^{-2}(0)/2\tau]\,\exp[-x^2/2D(0)\tau]\ . \tag{3.8}$$

Here the initial point is zero, and $D$ is the respective diagonal element of $\underline{D}$. The *mode* equals zero, while the mean (3.5) results by the asymmetric term with $x^3$.

Suppose that a given interval $[0, T]$ is partitioned into intervals of a length $\tau$, and consider the sum of the increments (3.7) . Stochastic calculus states that the limit $\tau \to 0$ is unchanged when $(\Delta W_i)^2$ is replaced by its expectation $\tau$; the sum thus converges to the anti-Ito integral.
*Comment*: While the simpler increments $\Delta X^i = \sqrt{D^{ii}}\, \Delta W_i + (D^{ii}{}_{,i}/2)\, \tau + o(\tau)$ lead to the same integral, they are Gaussian distributed and *disagree* with (3.8) , i.e. with the FPE.

Physically the case $\alpha = 1$ applies for tagged molecules in a fluid (with a temperature gradient): then $p$ is proportional to their density, and in view of (A.3) the probability (i.e. particle) current $\vec{J}$ obeys Fick's law, in each component when $\underline{D}$ is diagonal (a diagonal $\underline{D}$ also excludes a rotation of $\vec{J}$ around a maximum point of $p$, see (A.8) and the respective comment there). This specifies a local equilibrium, despite the temperature gradient.
In contrast to the linear case, the interaction with the other molecules is fully reciprocal; the state of the other molecules is given by the local equilibrium which is sufficient for the Markov property. Note that (3.8) exhibits a trend to the *warmer* side.

A nonzero $\vec{a}$ contributes $a^i(\vec{x})\tau$ to $\Delta X^i$, and the FPE reads

$$p_{,t} = \nabla \cdot [-\vec{a}p + (1/2)\, \underline{D}\, \nabla p] \tag{3.9}$$

(with the diagonal $\underline{D}$); the above densities are shifted by $a^i\tau$ (mind again the exception at an equilibrium point of $\vec{a}$).

## IV. *Supplement:* On steady states and weak noise

Steady states are possible solutions of a FPE with $p_{,t} \equiv 0$. They are associated with a steady current $\vec{J}_s$ , given by

$$\vec{J}_s = \vec{a}p - (1/2)\, \underline{D}\, \nabla p \tag{4.1}$$

for the quadratic case; in the linear case the righthand side is completed by $-(D^{ik}{}_{,k}/2)\, p$.

Clearly, $\vec{J}_s$ must satisfy $\nabla \cdot \vec{J}_s = 0$ (the FPE), as well as $\vec{J}_s \cdot \nabla p = 0$.

For ‘weak noise’ the extra term is neglected [13]. The steady states are thus formally the same in both cases, but the neglection of $-(D^{ik}{}_{,k}/2)\, p$ actually amounts to a constant $\underline{D}$.

## V. Summary and discussion

A differential approach (with generalized derivatives) saves the usual integral form; stochastic integrals are only mentioned to show the relation with existing work. Instead, functions of the Wiener process are considered, which lead to Markov chains for computing the resolving paths. Multiplicative noise entails two qualitatively different solutions that correspond to ‘Ito’ and to ‘anti-Ito’ (while $0 < \alpha < 1$ excludes the Markov property). Applied to particles in a medium with a temperature gradient, the first one holds for Brownian particles shifted to the cooler side (by stronger hits from the warm side), while the second one describes the motion of a ‘tagged’ molecule, moving to the warmer side (the latter case is associated with a local equilibrium). Generally, the criterion for choosing $\alpha$ is a possible feedback to the noise source: none for $\alpha = 0$ and maximum for $\alpha = 1$.

## *Appendix*: **A Fokker-Planck equation and its basic solution**

Consider the FPE (without $\vec{a}$, for simplicity)

$$p_{,t} = (1/2)\, \nabla \cdot (\, \underline{D}\, \nabla p) \tag{A.1}$$

with any symmetric and nonnegative $\underline{D}(\vec{x})$. The mean of an increment in $[0, \tau]$ is given by

$$< \Delta X^i > = (D^{ik}{}_{,k}/2)\, \tau + o(\tau)\ ; \tag{A.2}$$

this follows by the more usual form $p_{,t} = (1/2)[-D^{ik}{}_{,k} p + (D^{ik} p)_{,k}]_{,i}$ which equals (A.1) in view of $(D^{ik} p)_{,k} = D^{ik}{}_{,k} p + D^{ik} p_{,k}$ . It holds for the anti-Ito case.

With the *probability current*

$$\vec{J} := -(1/2)\,\underline{D}\,\nabla p \tag{A.3}$$

(A.1) becomes $p_{,t} + \nabla\cdot\vec{J} = 0$ which expresses the conservation of probability, i.e. the continuity of the paths. Note that $\vec{J} = \vec{0}$ where $\nabla p = \vec{0}$.

In the Ito case (with $< \Delta X^i > = 0$) (A.1) becomes

$$p_{,t} = (1/2)\,(D^{ik}p)_{,ki}\;. \tag{A.4}$$

Here the probability current reads

$$J^i = -(D^{ik}{}_{,k}/2)\,p \quad \text{where } \nabla p = \vec{0}\,. \tag{A.5}$$

It is interesting to consider the solution of (A.1) , with an initial deltafunction at $\vec{x} = \vec{0}$, say, and for small $\tau$. A suitable ansatz is

$$p(\vec{x},\tau) = [2\pi\, det\underline{D}(\vec{0})\,\tau]^{-r/2}\exp[-\Phi(\vec{x})\,/2\tau] \tag{A.6}$$

where $r$ is the rank of $\underline{D}$; (A.6) applies in the (sub)space where $\underline{D}$ is nonsingular. The $\vec{x}$-dependence is confined to $\Phi$ which only appears in the exponent; $\Phi$ can be seen as a kind of a quasipotential. Inserting this into (A.1) , dividing by $p$, and retaining only the leading terms for small $\tau$ (proportional to $\tau^{-2}$) leads to the equation for $\Phi$ :

$$\Phi = (1/4)\,(\nabla\Phi)^T\underline{D}\,\nabla\Phi\,. \tag{A.7}$$

This is a covariant equation (in the sense of tensor analysis), since $\Phi$ is a scalar and $\underline{D}$ is a twice contravariant tensor [9]. It thus preserves its form in new variables $\vec{z}(\vec{x})$, in particular when $\vec{z}$ is given by $d\vec{z} = [\underline{D}(\vec{x})]^{-1/2}d\vec{x}$, which replaces $\underline{D}$ by the unit matrix $\underline{I}$ : then $\Phi(\vec{z}) = (\vec{z})^T\underline{I}\,\vec{z}$ . The inverse transform yields $\Phi(\vec{x}) = (\vec{x})^T\underline{D}^{-1}\,\vec{x}$ , since $\underline{D}^{-1}$ is a twice covariant tensor which is also unity in $\vec{z}$. It follows that

$$p(\vec{x},\tau) = [2\pi\, det\underline{D}(\vec{0})\,\tau]^{-r/2}\exp[-\vec{x}^T\underline{D}^{-1}(\vec{x})\,\vec{x}/2\tau] \tag{A.8}$$

solves the FPE for small enough $\tau$. The mode is the starting point $\vec{x} = \vec{0}$, and the local curvature is not affected by an $\vec{x}$-dependence of $\underline{D}$. An asymmetry only arises by a cubic $x$- dependence of the exponent, and this determines the mean, see below. A nonconstant and nondiagonal $\underline{D}$ implies a rotational motion, since the evolution does then not follow the eigenvectors of $\underline{D}(\vec{0})$.

In one dimension one can infer from (A.8) that

$$p(x,\tau) \approx [2\pi D(0)\tau]^{-1/2}\,[1 + x^3 D'(0) D^{-2}(0)/2\tau]\,\exp[-x^2/2D(0)\tau] \qquad \text{(A.9)}$$

(by linearizing $D(x)$ in the exponent and by $\exp(x^3) \approx 1 + x^3$). This is normalized, and the term with $x^3$ yields the correct mean $[D'(0)/2]\,\tau$. The variance is $D(0)\tau + o(\tau)$.